\documentclass[review,3p]{elsarticle}

\usepackage{amsmath}
\usepackage{amsfonts}

\usepackage{hyperref}

\journal{Journal of Mathematical Analysis and Applications}

\DeclareMathOperator{\slope}{\mathtt{Slope}}  
\DeclareMathOperator{\Real}{Re}								
\DeclareMathOperator{\Imag}{Im}								

\newtheorem{theorem}{Theorem}
\newtheorem{lemma}{Lemma}

\newdefinition{remark}{Remark}
\newdefinition{claim}{Claim}

\newproof{proof}{Proof}

\begin{document}

\begin{figure}[ht]
  \includegraphics[scale=0.6]{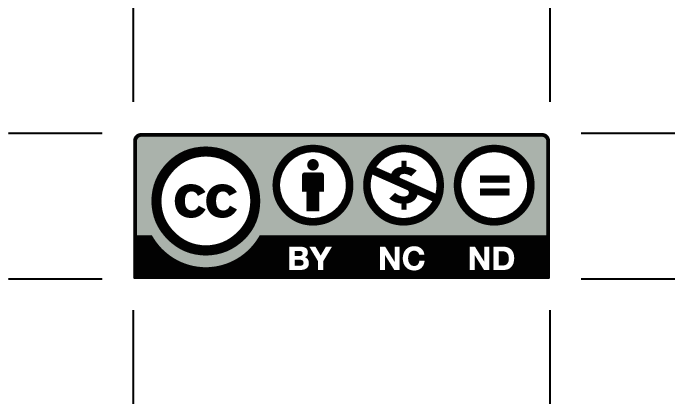}
\end{figure}

\begin{frontmatter}
  
\title{Trajectories of semigroups of holomorphic functions and harmonic measure\tnoteref{t1}}
\tnotetext[t1]{\textcopyright ~ 2019. This manuscript version is made available under the CC-BY-NC-ND 4.0 license \url{http://creativecommons.org/licenses/by-nc-nd/4.0/}}
\author{Georgios Kelgiannis\corref{mycorrespondingauthor}}
\cortext[mycorrespondingauthor]{Corresponding author}
\address{Department of Mathematics, Aristotle University of Thessaloniki, 54124, Thessaloniki, Greece}
\ead{gkelgian@math.auth.gr}
\ead[url]{https://users.auth.gr/gkelgian}
 
\begin{abstract}
  We give an explicit relation between the slope of the trajectory of a semigroup of holomorphic functions and the harmonic measure of the associated planar domain $ \varOmega $. We use this to construct a semigroup whose slope is an arbitrary interval in $ [ -\pi / 2,~ \pi / 2 ] $. The same method is used for the slope of a backward trajectory approaching a super-repulsive fixed point.
\end{abstract} 

  
\end{frontmatter}


\section{Semigroups of Holomorphic Functions}

A one-parameter continuous semigroup of holomorphic self-mappings of the unit disk $ \mathbb{D} $ is a family $ (\phi_t)_{t \in [0,\infty)},  $ such that:
\begin{enumerate}[(i)]
  \item $ \phi_{t+s} = \phi_t \circ \phi_s $, for all $ t,s \in [0,+\infty) $
  \item $ \phi_0 (z) = z $ 
  \item $  \displaystyle \lim_{t \rightarrow s} \phi_t (z) = \phi_s (z) $, for all $ s \in [0,+\infty) $.
\end{enumerate}

We will simply call $ (\phi_t) $ a semigroup. For general reference on semigroups we point to \cite{MR1098711}, \cite{MR2683159} and \cite{MR1849612}.

A semigroup is called \textit{elliptic} if it is not a group of hyperbolic rotations and it has an interior fixed point, which must be the same for all $ \phi_t,~ t > 0 $. If $ (\phi_t) $ is a non-elliptic semigroup, then there exists a unique 
point $ \xi \in \partial \mathbb{D} $, called the Denjoy-Wolff point of the semigroup \cite{MR0480965}, such that
\begin{equation} \label{DWlim}
  \lim_{t \rightarrow \infty} \phi_t(z) = \xi,~\text{ for every } z \in \mathbb{D}.
\end{equation}

A semigroup with no interior fixed point is called \textit{non-elliptic}. From now on we will only deal with non-elliptic semigroups. An important tool in the study of non-elliptic semigroups is the corresponding Koenigs function, see \cite{MR1098711}, \cite{MR2683159}, \cite{MR1849612} and the references therein. To every non-elliptic semigroup $ (\phi_t) $, corresponds a conformal mapping $ h : \mathbb{D} \rightarrow \varOmega $ such that:
\begin{enumerate}[(i)]
  \item $ h(\mathbb{D}) = \varOmega $, and
  \item $ h( \phi_t (z) ) = h(z) + t,~ z \in \mathbb{D},~ t \ge 0 $.
\end{enumerate}

The domain $ \varOmega $ is called the associated planar domain of $ (\phi_t) $. A domain $ \varOmega $ is called convex in the positive direction when $ \{ z + t : z \in \varOmega \} \subset \varOmega $, for all $ t \in [0,\infty) $. Obviously the associated planar domain of a semigroup is convex in the positive direction. The converse is also true; for every simply connected domain $ \varOmega $ convex in the positive direction, define
\begin{equation*}
  \phi_t(z) = h^{-1}(h(z)+t),
\end{equation*}
where $ h $ is the Riemann map that maps $ \mathbb{D} $ onto $ \varOmega $. It is easy to verify that $ (\phi_t) $, as defined above, is a semigroup. 

We are interested in the boundary fixed points of $ \phi_t $. These are defined using the notion of angular limit. When $ \phi(z) \rightarrow w' $ as $ z \rightarrow w $ through any sector at $ w $ we say that $ w' $ is the angular limit of $ \phi $ as $ z $ tends to $ w $; we write
\begin{equation*}
  \angle \lim_{z \rightarrow w} \phi(z) = w'.
\end{equation*}
A point $ w \in \partial \mathbb{D} $ is called a boundary fixed point of $ \phi $, when $ \angle \lim_{z \rightarrow w} \phi(z) = w $. For a boundary fixed point $ w $, we define the angular derivative at $ w $ to be
\begin{equation*}
  \phi'(w) = \angle \lim_{z \rightarrow w} \frac{w-\phi(z)}{w-z}.
\end{equation*}
In the case when $ \phi(\mathbb{D}) \subset \mathbb{D} $, we know \cite[p.82]{MR1217706} that $ \phi'(w) $ always exists and belongs to $ (0,+\infty) \cup \{ \infty \} $. Boundary fixed points in this case are divided into three categories; see \cite{MR2225072} and references therein.
\begin{enumerate}[(i)]
  \item When $ \phi'(w) \in (0,1] $, $ w $ is called an attractive point,
  \item when $ \phi'(w) \in (1,+\infty) $, $ w $ is called a repulsive point and
  \item when $ \phi'(w) = \infty $, $ w $ is called a super-repulsive point.
\end{enumerate}

The Denjoy-Wolff Theorem guarantees that, in the context of semigroups, the Denjoy-Wolff point $ \xi $ in relation (\ref{DWlim}), is the unique attractive boundary fixed point of $ \phi_t $, for all $ t > 0 $.
  
Non-elliptic semigroups can be categorized according to properties of the associated planar domain $ \varOmega $; see e.g. \cite{MR3451236}. Namely:
\begin{enumerate}[(i)]
  \item When $ \varOmega $ is contained in a horizontal strip, the semigroup is called hyperbolic.
  \item When $ \varOmega $ is not contained in a horizontal strip, but it is contained in a horizontal half-plane, the semigroup is called parabolic of positive hyperbolic step.
  \item When $ \varOmega $ is not contained in any horizontal half-plane, the semigroup is called parabolic of zero hyperbolic step.
\end{enumerate}

The trajectory of $ z \in \mathbb{D} $ of a semigroup $ ( \phi_t ) $ is defined as the curve
\begin{equation*}
  \gamma_{z} : [0,+\infty) \rightarrow \mathbb{D}, ~ \gamma_{z}(t) = \phi_t(z).
\end{equation*}

By utilizing the associated domain $ \varOmega $, every trajectory can be extended as follows. Let $ T $ be the infinum of $ \{ t : h(z) + t \in \varOmega \} $.  The extended trajectory of $ z $ is the curve defined by
\begin{equation}
  \gamma_{z} : (T,+\infty) \rightarrow \mathbb{D}, ~ \gamma_{z}(t) = h^{-1}( h(z) + t ).
\end{equation}
From now on $ \gamma_{z} $ will be used for the extended trajectory. In accordance with \cite{MR2225072}, we will define the $ \alpha $ and $ \omega $ limits of curves. For every curve $ \Gamma : (s_1,s_2) \rightarrow \mathbb{C} $, if there exists a strictly increasing sequence $ t_n  \rightarrow s_2 $, such that $ \Gamma(t_n) \rightarrow \xi $, then $ \xi $ is called an $ \omega $-limit point of $ \Gamma $. The set of all $ \omega $-limit points of $ \Gamma $ is called the $ \omega $-limit set and denoted by $ \omega(\Gamma) $. Replacing $ s_2 $ with $ s_1 $ and considering strictly decreasing sequences, we similarly define the $ \alpha $-limit point and the $ \alpha $-limit set $ \alpha(\Gamma) $.
From (\ref{DWlim}) it is obvious that for all $ z \in \mathbb{D} $ we have $ \omega(\gamma_{z}) =\{ \xi \} $, where $ \xi $ is the Denjoy-Wolff point. The set $ \alpha(\gamma_{z}) $ is also a single point which can be one of the following \cite{MR2225072}:
\begin{enumerate}[(i)]
  \item The point in $ \partial \mathbb{D} $ that corresponds to $ h(z)+T \in \partial \varOmega $, when $ T > -\infty $. 
  \item A boundary fixed point of $ (\phi_t) $, including the Denjoy-Wolff point $ \xi $, when $ T = -\infty $.
\end{enumerate}
An interesting problem is the study of the slope of $ \gamma_{z} $ as it approaches the boundary of $ \mathbb{D} $. For every $ \gamma_{z} $, we consider the corresponding curve
\begin{equation}
  t \in (T,+\infty) \rightarrow \arg(1-\bar{\xi}\gamma_{z}(t)) \in (-\frac{\pi}{2},\frac{\pi}{2}).
\end{equation}
The $ \omega $-limit set of the above curve will be the set of slopes of $ \gamma_{z} $ as it approaches the Denjoy-Wolff point $ \xi $ and it will be denoted by $ \slope^+ ( \gamma_{z} ) $. 
If $ \alpha(\gamma_{z}) = \{ \chi \} $ then similarly consider the curve
\begin{equation}
  t \in (T,+\infty) \rightarrow \arg(1-\bar{\chi}\gamma_{z}(t)) \in (-\frac{\pi}{2},\frac{\pi}{2}).
\end{equation}
The $ \alpha $-limit set of the above curve will be called the set of slopes of the backward trajectory $ \gamma_{z} $ as it approaches the boundary point $ \chi $ and it will be denoted by $ \slope^-( \gamma_{z} ) $. The following is already known about the $ \slope^+ ( \gamma_{z} ) $.
\begin{enumerate}[(i)]
  \item When a semigroup is hyperbolic, $ \slope^+ ( \gamma_{z} ) $ is a singleton depending on $ z $. 
  \item When a semigroup is parabolic of positive hyperbolic step, $ \slope^+ ( \gamma_{z} ) $ is either $ \{ \pi / 2 \} $ or $ \{ - \pi / 2 \} $ and it is independent of $ z $.
\end{enumerate}
When a semigroup is parabolic of zero hyperbolic step, it was conjectured that 
$ \slope^+ ( \gamma_{z} ) $ is again a singleton. This was proven but only under some additional assumptions, see e.g.\ \cite{MR3240989} and \cite{MR2731701}.
The existence of a semigroup with $ \slope^+ ( \gamma_{z} ) = [-\pi/2,\pi/2] $ was first proven in \cite{MR3441527} and \cite{MR3318309}. In a more recent result, Bracci et al. \cite{bracci2018non} show that there exists a semigroup such that $ \slope^+ ( \gamma_{z} ) \subset (-\pi/2,\pi/2) $ but it is not a singleton. Also in \cite{bracci2018asy} we find an example with $ \slope^+ ( \gamma_{z} ) = [-\pi/2,\alpha] $, for some $ -\pi/2 < \alpha < \pi / 2 $. 

In \cite{MR3318309} the authors posed the problem of constructing examples of one-parameter semigroups $ (\phi_t) $ with $ \slope^+ (\gamma_{z}) = [\theta_1, \theta_2] $ for any given $ \theta_1, \theta_2 \in [-\pi/2, \pi/2 ],~ \theta_1 < \theta_2 $. We will construct such a semigroup.
\begin{theorem}\label{slope+}
  If $ \theta_1 < \theta_2 $ are real numbers with $ |\theta_j| \le \pi/2,~ j=1,2 $, then 
  there exists a semigroup of holomorphic functions $ (\phi_t) $ such that
  \begin{equation}
    \slope^+ (\gamma_{z}) = [\theta_1,\theta_2].
  \end{equation} 
\end{theorem} 

For the $ \slope^- ( \gamma_{z} ) $ similar results were only known for the following cases \cite{MR2225072}:
\begin{enumerate}[(i)]  
  \item When the  $ \alpha $-limit of $ \gamma_{z} $ is the Denjoy-Wolff point  $ \xi $, $ \slope^- ( \gamma_{z} ) $ is a singleton, which is either $ \{ \pi/2 \} $ or $ \{ -\pi/2 \} $.
  \item When the  $ \alpha $-limit of $ \gamma_{z} $ is a repulsive point, $ \slope^- ( \gamma_{z} ) $ is a single point, which belongs in $ ( -\pi/2 , \pi/2 ) $.  
\end{enumerate}

We prove that, in the case of super-repulsive points, a semigroup can have a wildly oscillating trajectory, quite similar to the case of a 
parabolic semigroup of zero hyperbolic step.
\begin{theorem}\label{slope-}
  If $ \theta_1 \le \theta_2 $ are real numbers with $ |\theta_j| \le \pi/2,~ j=1,2 $, then there exists a semigroup of holomorphic functions $ (\phi_t) $ and a point $ z \in \mathbb{D} $, such that the $ \alpha $-limit of $ \gamma_{z} $ is a super-repulsive point and
  \begin{equation*}
    \slope^- (\gamma_{z}) = [\theta_1,\theta_2].
  \end{equation*}
\end{theorem}

\section{Harmonic measure}

To prove the aforementioned results we need to establish a relationship between the slope of a trajectory $ \gamma_ {z} $ and certain harmonic measures in the associated planar domain $ \varOmega $ of a semigroup.

The harmonic measure is the solution $ u $ of the generalized Dirichlet problem for the Laplacian in a domain $ D $, with boundary values equal to $ 1 $ on $ E \subset \partial \varOmega $ and $ 0 $ on $ \partial \varOmega \setminus E $. We will be using the notation $  \omega(z,E,D) $.

Two basic properties of the harmonic measure that we will use are conformal invariance and domain monotonicity. When $ \phi: \mathbb{D} \rightarrow \varOmega $ is a conformal map, we know that, if $ A $ is the set of accessible points of $ \partial \varOmega $, we can extend $ \phi^{-1} $ to $ A $. In that sense, when $ E \subset A $ we have \cite[p.206]{MR2450237}
\begin{equation}
  \omega(z,\phi^{-1}(E),\mathbb{D}) = \omega( \phi(z), E,\varOmega ).
\end{equation}
This implies that when an arc $ \widehat{ab} \subset \partial \mathbb{D} $ corresponds, through $ \phi $, to a boundary set $ E \subset \partial \varOmega $, in the sense of Caratheodory boundary correspondence, then 
\begin{equation}\label{confInv}
  \omega(z,\widehat{ab},\mathbb{D}) = \omega( \phi(z), E, \varOmega ).
\end{equation}
When for two domains $ D_1,D_2 $ in $ \mathbb{C}_\infty $, with $ D_1 \subset D_2 $, we have a set $ B \subset \partial D_1 \cap \partial D_2 $, then \cite[p. 102]{MR1334766}
\begin{equation}
  \omega(z,B,D_1) \le \omega(z,B,D_2).
\end{equation} 
We also know that \cite[p.155]{MR0060009}, if $ \widehat{ab} \subset \partial \mathbb{D} $ is a circular arc, then the level set 
\begin{equation}
  L_k = \{ \zeta \in \mathbb{D} : \omega(\zeta,\widehat{ab},\mathbb{D}) = k \},~0<k<1,
\end{equation}
is a circular arc with endpoints $ a $ and $ b $ that meets the unit circle with angle $ k\pi $. We will also use the notation
\begin{equation}
  \widehat{L}_k = \{ \zeta \in \mathbb{D} : \omega(\zeta,\widehat{ab},\mathbb{D}) > k \}. 
\end{equation}

In order to establish a relation between certain harmonic measures in the case when $ D $ contains, in a specific way, a rectangle, we introduce the following notation.

For any set $ B $ in the complex plane $ \mathbb{C} $, let $ B^+ = B \cap \{ z: \Imag z > 0 \} $ and 
$ B^- = B \cap \{ z: \Imag z < 0 \} $. Let
\begin{equation} 
  S_d = \{z: -d< \Imag z < d \}
\end{equation}
be a horizontal strip of width $ 2d $,
\begin{equation}
  S_{d,u} = \{ z: -d< \Imag z < d,~ -u < \Real z < u \}
\end{equation}
be a rectangle centered at the origin with width $ 2d $ and length $ 2u $,
\begin{equation}
  B_{d,u} = \{ z: \Imag z = d,~ -u < \Real z < u \} 
\end{equation}
be the upper side of $ S_{d,u} $ and $ B_{-d,u} $ be the lower side of $ S_{d,u} $. Betsakos \cite{MR3441527} has proven the following:
\begin{lemma}\label{origL}
  Let $ \varOmega $ be a planar domain, convex in the positive direction. Assume that $ \mathbb{R} \subset \varOmega $ and that $ (\partial \varOmega)^+ \ne \emptyset,~ (\partial \varOmega)^- \ne \emptyset $. Let $ \epsilon > 0 $ and $ d > 0 $. There exists a $ u_0 > 0 $ with the property: 
  If $ y \in (-d,d), S_{d,u_0} \subset \varOmega $ and $ B_{d,u_0} \cup B_{-d,u_0} \subset \partial \varOmega $, then
  \begin{equation}\label{origR}
    |\omega(iy,(\partial \varOmega)^+, \varOmega) - \omega(iy,(\partial S_d)^+,  S_d) | < \epsilon.
  \end{equation}
\end{lemma}

In the original proof $ \varOmega $ is fixed. However, a close inspection of the proof shows that $ u_0 $ depends only on $ d $, not on the set $ \varOmega $ and that (\ref{origR}) holds for all $ u > u_0 $. We will use a variation of Lemma \ref{origL}.

For $ w \in \mathbb{C} $, $ d_1,d_2,u > 0 $, we consider the rectangles 
\begin{equation*}
  A(w,d_1,d_2,u) = \{ x+iy : |x - \Real w| < u/2 ,~ \Imag w - d_2 < y < \Imag w + d_1 \}.
\end{equation*}  
Let also, for $ A = A(w,d_1,d_2,u) $,
\begin{equation*}
  \partial_h A = \{ x+iy :|x - \Real w| < u/2 ,~ y = \Imag w - d_2 \text{ or } y = \Imag w + d_1 \},
\end{equation*}
be the horizontal border of $ A $. Finally for $ z \in \mathbb{C} $, let
\begin{equation} 
  \partial_{z}^+\varOmega = \partial \varOmega \cap \{ \zeta : \Imag \zeta > \Imag z \} 
\end{equation}
be the part of the border of $ \varOmega $ that lies above $ z $. Note that when $ z \in \mathbb{R} $ we have $ \partial_{z}^+\varOmega =(\partial \varOmega)^+ $.
Note also that if the distances of $ iy $ from the upper and lower parts of a strip are respectively $ d_1 $ and $ d_2 $, by applying standard conformal maps, one can see that
\begin{equation}
  \omega(iy,(\partial S_d)^+, S_d) = \dfrac{d_2}{d_1+d_2}.
\end{equation} 
By conformal invariance of the harmonic measure, Lemma \ref{origL} can be restated as follows.
\begin{lemma}\label{stripap}
  Let $ d_1,d_2 > 0 $. Then for every $ \epsilon > 0 $, there exists a $ u_0 > 0 $, such that for every $ u > u_0 $ and for all domains $ \varOmega $, convex in the positive direction, the following property holds:
  If $ A = A(w,d_1,d_2,u) \subset \varOmega $ and $ \partial_h A \subset \partial \varOmega $, then 
  \begin{equation}
    \left| \omega(w,\partial_{w}^+\varOmega, \varOmega) - \frac{d_2}{d_1+d_2} \right| < \epsilon.
  \end{equation}
\end{lemma}

We will be working with domains convex in the positive direction but we point out that by a small modification of the proof found in \cite{MR3441527}, we can drop this requirement.

Let $ z \in \mathbb{D} $. We will prove that the slope of the trajectory  $ \gamma_{z} $ of a semigroup of holomorphic functions $ (\phi_t) $ is determined by certain harmonic measures. Consider the function  
\begin{equation}\label{omt}
  \omega_z(t) = \omega(h(z)+t,\partial_{h(z)}^+\varOmega,\varOmega),~t\in(0,+\infty).
\end{equation}
Betsakos \cite{MR3441527} constructed a semigroup such that for every $ z \in \mathbb{D},~ \slope^+(\gamma_{z}) = [-\pi/2,\pi/2] $, by considering the behavior of $ \omega_0(t) $ as $ t \rightarrow + \infty $. We will prove an explicit relation between the behavior of  $ \omega_z(t) $ and the slopes of $ \gamma_{z} $. We will then use it to construct a semigroup such that $ \slope^+(\gamma_{z}) = [\theta_1,\theta_2] $ with $ -\pi/2 \le \theta_1 < \theta_2 \le \pi/2 $. The same principles will be extended to an analogous result for the $ \slope^-(\gamma_{z}) $.
\begin{theorem}\label{harmonic-slope}
  Let $ (\phi_t) $ be a semigroup of holomorphic functions in $ \mathbb{D} $. Denote by $ h $ the corresponding Koenigs function and by $ \varOmega = h(\mathbb{D}) $ the associated planar domain. For $ z \in \mathbb{D} $, with $ \partial_{h(z)}^+\varOmega \ne \emptyset $ and $ \partial_{h(z)}^+\varOmega \ne \partial \varOmega $, let $ a_1 = \limsup_{t \rightarrow \infty} \omega_z(t) $ and $ a_2 = \liminf_{t \rightarrow \infty} \omega_z(t) $. Then
  \begin{equation}\label{g+}
    \slope^+ (\gamma_{z}) = \left[ \pi(1/2-a_1) , ~  \pi(1/2-a_2) \right].
  \end{equation}	
  If, in addition, for that $ z $, the trajectory $ \gamma_{z} $ is defined for all $ t \in (-\infty, 0] $ and we have $ b_1 = \limsup_{t \rightarrow -\infty} \omega_z(t) $ and $ b_2 = \liminf_{t \rightarrow -\infty} \omega_z(t) $, then
  \begin{equation}\label{g-}
    \slope^- (\gamma_{z}) = \left[ \pi(1/2-b_1) , ~  \pi(1/2-b_2) \right].
  \end{equation}
\end{theorem}

Using the above theorem we can argue about the slopes of the trajectories of $ (\phi_t) $ by focusing on the image $ h(\mathbb{D}) $ and looking at the behavior of the harmonic measure  on the points of the half-line $ \{ h(z) + t : t > 0 \} $, or on $ \{ h(z) - t : t > 0 \} $ for the backward trajectories.

\section{Proofs}

\begin{proof}[Theorem \ref{harmonic-slope}]
 
We assume that the Denjoy-Wolff point of $ (\phi_t) $ is $ \xi $ and the $ \alpha$-limit of $ \gamma_{z} $ is $ \chi $. 
Let $ \widehat{\chi \xi} $ be the arc on $ \partial \mathbb{D} $ between $ \chi $ and $ \xi $, corresponding through $ h(z) $ to $ \partial_{h(z)}^+\varOmega $. Note that $ \partial_{h(z)}^+\varOmega \ne \emptyset $ and $ \partial_{h(z)}^+\varOmega \ne \partial \varOmega $ imply $ \chi \ne \xi $. Also since $ h $ is conformal we have that $ \widehat{\chi \xi} $ is the arc that runs clockwise from $ \chi $ to $ \xi $. We know that the level set 
\begin{equation}
  L_k = \{ \zeta \in \mathbb{D} : \omega(\zeta,\widehat{\chi \xi},\mathbb{D}) = k \},~0<k<1,
\end{equation}
is a circular arc with endpoints $ \chi $ and $ \xi $ that meets the unit circle with angle $ k\pi $. 
Let $ \widehat{L}_k = \{ \zeta \in \mathbb{D} : \omega(\zeta,\widehat{\chi \xi},\mathbb{D}) > k \} $ and $ \Gamma_k $ be the half-line emanating from $ \xi $ that is tangent to $ L_k $ at $ \xi $. 
If $ \zeta $ lies on $ \Gamma_k $ then $ \arg(1-\overline{\xi}\zeta) =  \pi/2 - \pi k  = \pi(1/2-k) $.

By conformal invariance of the harmonic measure (\ref{confInv}), 
\begin{equation}\label{pullback}
  \omega_z(t) = \omega(h(z) + t,\partial_{h(z)}^+\varOmega,\varOmega) = \omega(\phi_t(z),\widehat{\chi \xi},\mathbb{D}). 
\end{equation}
Let $ a_1 = \limsup_{t \rightarrow \infty} \omega_z(t) $ and $ \theta_1 = \pi(1/2-a_1) $ the corresponding angle.

We will prove that $ \theta_1 = \min\{\slope^+(\gamma_{z})\} $.

\begin{claim} If $ \theta \in \slope^+(\gamma_{z}) $ then $ \theta_1 \le \theta $.

If $ a_1 = 1 $ then $ \theta_1 = -\pi/2 $ and we are done. If not, 
since $ a_1 = \limsup_{t \rightarrow \infty} \omega_z(t) $, from (\ref{pullback}) we must also have
\begin{equation}
  \limsup_{t \rightarrow \infty} \omega(\phi_t(z),\widehat{\chi \xi},\mathbb{D}) = a_1.
\end{equation}
Assume that $ \theta \in \slope(\gamma_{z}) $ with $ \theta_1 > \theta = \pi(1/2-a) $.
So there is an $ \epsilon > 0 $ such that $ a_1 < a_1 + \epsilon / 2 < a_1 + \epsilon < a $.
Then there is a sequence $ t_n \rightarrow \infty $ such that all but finite of the points $ \phi_{t_n}(z) $ lie above $ \Gamma_{a_1 + \epsilon} $ for some  $ \epsilon>0 $. This means that $ \phi_{t_n}(z) \in \widehat{L}_{a_1+ \epsilon/2} $ for almost all $ n $. This implies that $ \lim_{t_n\rightarrow\infty}\omega(\phi_{t_n}(z),\widehat{\chi \xi},\mathbb{D}) \ge a_1+ \epsilon/2 $, a contradiction. So $ \theta_1 \le \theta $.    
\end{claim}

\begin{claim} $ \theta_1 \in \slope^+(\gamma_{z}) $. 

Since there exists $ t_n $ with $ \omega(\phi_{t_n}(z),\widehat{\chi \xi},\mathbb{D}) \rightarrow a_1 $ we have that  $  \arg(1-\overline{\xi}\phi_{t_n}(z)) \rightarrow \theta_1 $ and so $ \theta_1 \in \slope^+(\gamma_{z}) $.

\end{claim}

We have shown that $ \theta_1 = \min\{\slope^+(\gamma_{z})\} $. Using the same arguments we can show that if $ a_2 = \liminf_{t \rightarrow \infty} \omega_z(t) $ and $ \theta_2 = \pi(1/2-a_2) $ we have $ \theta_2 = \max\{\slope^+(\gamma_{z})\} $. This means that $ \slope^+(\gamma_{z}) = \left[ \pi(1/2-a_1) , ~  \pi(1/2-a_2) \right] $.

In the case when the $ \alpha $-limit of $ \gamma_{z} $ is a super-repulsive point, replacing $ \infty $ with $ -\infty $ and $ \xi $ with $ \chi $, using the same arguments, we obtain relation (\ref{g-}) for the $ \slope^- (\gamma_{z}) $.

\end{proof}

\begin{remark}
  The only property of the set $ \partial_{h(z)}^+\varOmega $ that we use is that it corresponds, through $ h^{-1} $, to an arc $  \widehat{\chi \xi} $  on $ \partial \mathbb{D} $ with $ \xi $ being the Denjoy-Wolff point, or $ \chi $ being the $\alpha$-limit of $ \gamma_{z} $, and $ \chi \ne \xi $. This means that even when $ \partial_{h(z)}^+\varOmega = \emptyset $ or $ \partial_{h(z)}^+\varOmega = \partial \varOmega $ we can use the same approach by choosing a suitable subset of $ \partial \varOmega $. 
\end{remark}

\begin{proof}[Theorem \ref{slope+}]
We will only prove the result for $ |\theta_1|, |\theta_2| < \pi/2 $ for simplicity. Small variations of the proof can also account for the cases of $ \theta_1 = -\pi / 2 $ or $ \theta_2 = \pi / 2 $. We will essentially present these variations in the proof of Theorem \ref{slope-}. We will modify the construction found  in \cite{MR3441527} and construct a set $ \varOmega $ such that for the associated semigroup we have $ \slope^+ (\gamma_{0}) = [\theta_1,\theta_2] $. Let $ E[\zeta] = \{ \zeta+t:t\le 0 \} $ be the half-line, parallel to the real axis, starting from $ \zeta $ and extending to the left. Let $ a_1 = \frac{1}{2} - \frac{\theta_1}{\pi} $ and $ a_2 = \frac{1}{2} - \frac{\theta_2}{\pi} $, so that $ 0 < a_2 < a_1 < 1 $. Let $ r_n, \rho_n $ be sequences such that
\begin{equation}\label{a1def}
  r_n = \frac{1-a_1}{a_1} \rho_n \\
\end{equation}
and
\begin{equation}\label{a2def}
  r_n = \frac{1-a_2}{a_2} \rho_{n-1}, ~ n \ge 2.
\end{equation}
Since $ a_1 > a_2 $, both $ r_n $ and $ \rho_n $ are increasing. Note that these depend only on the choice of $ a_1, a_2 $ and $ r_1 $. For example $ a_1 = \dfrac{3}{4} $, $ a_2 = \dfrac{1}{3} $ and $ r_1 = 6 $ gives
\begin{equation*}
  r_n = 6^n \text{ and } \rho_n = 3 \cdot 6^n.
\end{equation*}
It is easy to see that definitions (\ref{a1def}) and (\ref{a2def}) indeed give
\begin{equation}\label{a1a2}
  \frac{\rho_n}{\rho_n+r_n} = a_1 \text{ and } \frac{\rho_{n-1}}{\rho_{n-1}+r_n} = a_2.
\end{equation}  	
Note that for $ w = 0 $ we have $ \partial_{w}^+\varOmega = (\partial \varOmega)^+ $ and choose an increasing sequence $ u'_n $ from Lemma \ref{stripap}, such that the following hold:\\
When $ n= 2k-1 $, for  all $ \varOmega $ with $ A = A(w,r_k,\rho_k,u'_n) \subset \varOmega $ and $  \partial A_h \subset \partial \varOmega $,
\begin{equation}\label{rel1}
  | \omega(w,(\partial \varOmega)^+,\varOmega) - \frac{\rho_k}{\rho_k+r_k}| < \frac{1}{n}
\end{equation}
and for  all $ \varOmega $ with $ A = A(x,r_{k+1},\rho_k,u'_n) \subset \varOmega $ and $  \partial A_h \subset \partial \varOmega $,  	
\begin{equation}\label{rel2}
  | \omega(w,(\partial \varOmega)^+,\varOmega) - \frac{\rho_k}{\rho_k+r_{k+1}}| < \frac{1}{n}.
\end{equation}
When $ n= 2k $, for  all $ \varOmega $ with $ A = A(x,r_{k+1},\rho_{k},u'_n) \subset \varOmega $ and $  \partial A_h \subset \partial \varOmega $,
\begin{equation}\label{rel3}
  | \omega(w,(\partial \varOmega)^+,\varOmega) - \frac{\rho_{k}}{\rho_{k}+r_{k+1}}| < \frac{1}{n}
\end{equation}
and for  all $ \varOmega $ with $ A = A(x,r_{k+1},\rho_{k+1},u'_n) \subset \varOmega $ and $  \partial A_h \subset \partial \varOmega $,  	
\begin{equation}\label{rel4}
  | \omega(w,(\partial \varOmega)^+,\varOmega) - \frac{\rho_{k+1}}{\rho_{k+1}+r_{k+1}}| < \frac{1}{n}.
\end{equation}

Consider the partial sums $ u_n = \sum_{j=1}^n u'_j $ and set
\begin{equation}
  \varOmega = \mathbb{C} \setminus \bigcup_{n=1}^{\infty}(E[u_{2k-1}+ir_k] \cup E[u_{2k}-i\rho_k]).
\end{equation}

The way $ \varOmega $ was constructed we have that $ \varOmega $ is convex in the positive direction. We also have that, for $ n = 2k-1 $, for the rectangles   
$ A=A(x_n,r_k,\rho_{k},u'_{n}) $ we have $ A \subset \varOmega $ and $  \partial A_h \subset \partial \varOmega $, where $ x_n = (u_{n}+u_{n-1})/2 $. Obviously   $ x_n \rightarrow \infty $.  For $ n = 2k $ the same holds for $ A=A(x_n,r_{k+1},\rho_{k},u'_n) $.

So for $ n = 2k - 1 $, from relations (\ref{a1a2}) and (\ref{rel1}), we have,
\begin{equation} 
  | \omega(x_n,(\partial \varOmega)^+,\varOmega) -  a_1 | < \frac{1}{n}
\end{equation}
and for $ n = 2k $, from relations (\ref{a1a2}) and  (\ref{rel3}),
\begin{equation} 
  | \omega(x_n,(\partial \varOmega)^+,\varOmega) - a_2 | < \frac{1}{n}. 
\end{equation}

So we have found two sequences $ x_{2k-1} \in \mathbb{R} $ and $ x_{2k} \in \mathbb{R} $ with respective limits $ a_2 $ and $ a_1 $. That means
\begin{equation} 
  [a_2,a_1] \subset [\liminf_{t \rightarrow \infty} \omega_0(t),\limsup_{t \rightarrow \infty} \omega_0(t)] .
\end{equation}
We proceed to show the opposite inclusion.
Consider a pair $ x_{2k-1}, x_{2k} $ on the real line. Note that the rectangles $ A(x_{2k-1},r_k,\rho_{k},u'_{2k-1})$ and $ A(x_{2k},r_{k+1},\rho_{k},u'_{2k}) $   are both contained in $ \varOmega $. 	
\begin{figure}[ht]
 	\label{fig2} 
 	\caption{A part of the set $ \varOmega $}
\centering
\includegraphics{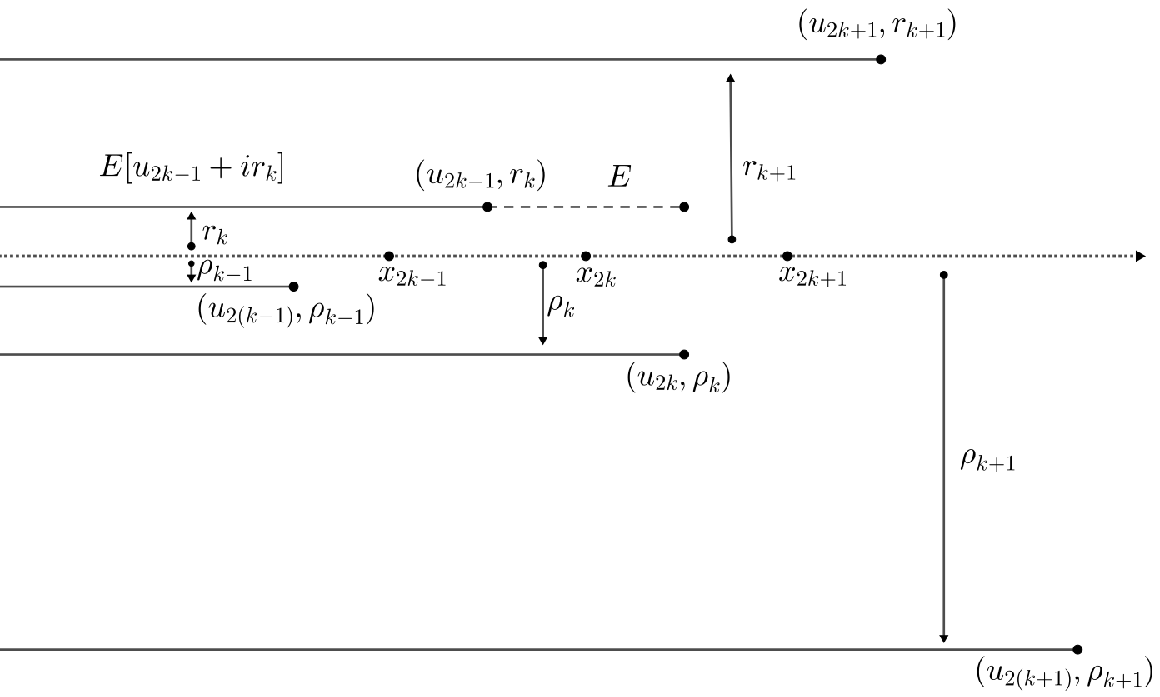}
\end{figure}

Consider the set $ \varOmega_1 = \varOmega \setminus E $, where $ E = \{ x+iy : y = r_k ,~ u_{2k-1} < x \le u_{2k} \} $. In Figure \ref{fig2}, $ E $ is the dotted segment. Obviously $ \varOmega_1 \subset \varOmega $ and $ (\partial \varOmega)^- = (\partial \varOmega_1)^- $. Also for all $ x \in [ x_{2k-1},x_{2k} ] $, since $ u'_n $ is increasing, we have that 
\begin{equation} 
  A=A(x,r_k,\rho_{k},u'_{2k-1}) \subset \varOmega_1 \text{ and } \partial A_h \subset \partial \varOmega_1. 
\end{equation}
Using the domain monotonicity of the harmonic measure and relation (\ref{rel1}) we get
\begin{align*}
  \omega(x,(\partial \varOmega)^+,\varOmega) = 1 - \omega(x,(\partial \varOmega)^-,\varOmega) &\le 1 - \omega(x,(\partial \varOmega_1)^-,\varOmega_1) \\
  &= \omega(x,(\partial \varOmega_1)^+,\varOmega_1) < a_1 + \frac{1}{n}. 
\end{align*}
Similarly consider $ \varOmega_2 = \varOmega \cup E[u_{2k-1}+ir_k] $. Again for all $ x \in [ x_{2k-1},x_{2k} ] $, we have $ A = A(x,r_{k+1},\rho_{k},u'_{2k-1}) \subset \varOmega_2 $ and $ \partial A_h \subset \partial \varOmega_2 $. Since $ \varOmega \subset \varOmega_2 $, considering (\ref{rel2}),
\begin{equation*}
  \omega(x,(\partial \varOmega)^+,\varOmega) > \omega(x,(\partial \varOmega_2)^+,\varOmega_2) > a_2 - \frac{1}{n}. 
\end{equation*}
We can likewise treat the case where $ x \in [ x_{2k},x_{2k+1} ] $.
These inequalities show that if there exists a sequence $ t_k \rightarrow \infty $ with $  \lim_{k \rightarrow \infty} \omega_0(t_k) = a $ then $ a_2 \le a \le a_1 $.

We have shown that $ a_1 =  \limsup_{t \rightarrow \infty}\omega_0(t) $ and $ a_2 = \liminf_{t \rightarrow \infty} \omega_0(t) $. Considering the semigroup $ (\phi_t) $ that corresponds to the set $ \varOmega $, the desired result follows from Theorem \ref{harmonic-slope}.

\end{proof}

\begin{proof}[Theorem \ref{slope-}]
As in the above proof let $ b_1 = \frac{1}{2} - \frac{\theta_1}{\pi} $, $ b_2 = \frac{1}{2} - \frac{\theta_2}{\pi} $ and $ r_n, \rho_n $ be sequences
such that
\begin{equation*}
  r_n = \frac{1-b_2}{b_2} \rho_n,
\end{equation*}
and
\begin{equation*}
  r_n = \frac{1-b_1}{b_1} \rho_{n-1}, ~ n \ge 2.
\end{equation*}
Since $ b_1 > b_2 $ we have that both $ r_n $ and $ \rho_n $ are decreasing sequences.
Note that these depend only on the choice of $ b_1,b_2 $ and $ r_1 $. Similar to the above proof, if for example $ b_1 = \dfrac{3}{4} $, $ b_2 = \dfrac{1}{3} $ and $ r_1 = \dfrac{1}{3} $, we get
\begin{equation*}
  r_n = \frac{1}{3} \cdot 6^{-(n-1)} \text{ and } \rho_n = 6^{-n}.
\end{equation*}

We define sequences $ u_n $, $ u_n' $ in the exact same way as in the proof of Theorem \ref{slope+}. This means that we can use relations (\ref{rel1} - \ref{rel4}).  Now $ \varOmega $ can be defined as
\begin{equation*}
  \varOmega = \mathbb{C} \setminus \bigcup_{n=1}^{\infty}(E[-u_{2k-1}+ir_k] \cup E[-u_{2k} + i\rho_k]).
\end{equation*}

Obviously $ \varOmega $ is convex in the positive direction and $ \gamma_{0} $ is defined for $ t \in (-\infty,+\infty) $. Similarly with before we take $ x_n = -(u_{n}+u_{n-1})/2 $. We have that $ x_n $ goes to $ -\infty $ and for the subsequences $ x_{2k-1} $ and $ x_{2k} $ we get
\begin{align*}
  \lim_{k \rightarrow \infty} \omega(x_{2k-1},(\partial \varOmega)^+,\varOmega) &= b_1 \text{ and } \\
  \lim_{k \rightarrow \infty} \omega(x_{2k},(\partial \varOmega)^+,\varOmega) &= b_2.
\end{align*}

We can show the opposite inclusion with the same arguments as in the proof of Theorem \ref{slope+}. Again from Theorem \ref{harmonic-slope} we get $	\slope^- (\gamma_{0}) = [\theta_1,\theta_2] $.

We will now consider the case when $ b_2 = 0 $. We modify our sequences so that
\begin{equation*}
  r_n =  (n+m) \rho_n,
\end{equation*}
and
\begin{equation*}
  r_n = \frac{1-b_1}{b_1}, ~ n \ge 2,
\end{equation*}	
where $ m $ is taken big enough, so that for all $ n $ we have $ n+m > \frac{1-b_2}{b_2} $. We again have two decreasing sequences. The proof works out in the same way except that now, for $ n=2k-1 $, relation (\ref{rel1}) becomes 
\begin{equation}
  \omega(x_n,(\partial \varOmega)^+,\varOmega) < \frac{1}{n+m+1} + \frac{1}{n} < \frac{2}{n}     
\end{equation}
for all $ n $. Obviously $ \omega(x_{2k-1},(\partial \varOmega)^+,\varOmega) \rightarrow 0 $ as $ k \rightarrow \infty $ and as before we have $ \omega(x_{2k},(\partial \varOmega)^+,\varOmega) \rightarrow b_1 $. 
Similarly in the case when $ b_1 = 1 $ we take
\begin{equation*}
  r_n = \frac{1-b_2}{b_2} \rho_n,
\end{equation*}
and
\begin{equation*}
  r_n = \frac{1}{n+m}, ~ n \ge 2,
\end{equation*}	
where $ m $ is taken big enough, so that for all $ n $ we have $ \dfrac{1}{n+m} < \dfrac{1-b_2}{b_2} $. As before, note that, for $ n=2k $, relation (\ref{rel3}) becomes 
\begin{align*}
  \omega(x_n,\partial \varOmega^+,\varOmega) &> \frac{n+m}{n+m+1 } - \frac{1}{n} \\
  &> \frac{n+m}{n+m+1 } - \frac{2}{n+m+1} = 1 - \frac{3}{n+m+1},
\end{align*}
for all $ n > m + 1 $. Obviously $ \omega(x_{2k},(\partial \varOmega)^+,\varOmega) \rightarrow 1 $ as $ k \rightarrow \infty $, while $ \omega(x_{2k-1},(\partial \varOmega)^+,\varOmega) \rightarrow b_2 $.

Combining the above we can also construct an example with $ \slope^- ( \gamma_{z} ) = [-\pi/2,\pi/2] $. Note that in this case we can simply use 
\begin{equation*}
  r_n =  n \rho_n,
\end{equation*}
and
\begin{equation*}
  r_n = \frac{1}{n}\rho_{n-1}, ~ n \ge 2,
\end{equation*}	
which coincides with what was used in \cite{MR3441527}.

\end{proof}

\section*{Acknowledgements}
  
I would like to thank professor D. Betsakos, my thesis advisor, for his help.
  
This research did not receive any specific grant from funding agencies in the public, commercial, or not-for-profit sectors.
  
Declarations of interest: none.
  
\bibliography{slope}

\begin{thebibliography}{16}
\expandafter\ifx\csname natexlab\endcsname\relax\def\natexlab#1{#1}\fi
\providecommand{\url}[1]{\texttt{#1}}
\providecommand{\href}[2]{#2}
\providecommand{\path}[1]{#1}
\providecommand{\DOIprefix}{doi:}
\providecommand{\ArXivprefix}{arXiv:}
\providecommand{\URLprefix}{URL: }
\providecommand{\Pubmedprefix}{pmid:}
\providecommand{\doi}[1]{\href{http://dx.doi.org/#1}{\path{#1}}}
\providecommand{\Pubmed}[1]{\href{pmid:#1}{\path{#1}}}
\providecommand{\bibinfo}[2]{#2}
\ifx\xfnm\relax \def\xfnm[#1]{\unskip,\space#1}\fi
\bibitem[{Abate(1989)}]{MR1098711}
\bibinfo{author}{Abate, M.} (\bibinfo{year}{1989}).
\newblock {\it \bibinfo{title}{Iteration theory of holomorphic maps on taut
  manifolds}\/}.
\newblock Research and Lecture Notes in Mathematics. Complex Analysis and
  Geometry.
\newblock \bibinfo{publisher}{Mediterranean Press, Rende}.
\bibitem[{Berkson \& Porta(1978)}]{MR0480965}
\bibinfo{author}{Berkson, E.}, \& \bibinfo{author}{Porta, H.}
  (\bibinfo{year}{1978}).
\newblock \bibinfo{title}{Semigroups of analytic functions and composition
  operators}.
\newblock {\it \bibinfo{journal}{Michigan Math. J.}\/},  {\it
  \bibinfo{volume}{25}\/}, \bibinfo{pages}{101--115}.
  \DOIprefix\doi{10.1307/mmj/1029002009}.
\bibitem[{Betsakos(2016{\natexlab{a}})}]{MR3451236}
\bibinfo{author}{Betsakos, D.} (\bibinfo{year}{2016}{\natexlab{a}}).
\newblock \bibinfo{title}{Geometric description of the classification of
  holomorphic semigroups}.
\newblock {\it \bibinfo{journal}{Proc. Amer. Math. Soc.}\/},  {\it
  \bibinfo{volume}{144}\/}, \bibinfo{pages}{1595--1604}.
  \DOIprefix\doi{10.1090/proc/12814}.
\bibitem[{Betsakos(2016{\natexlab{b}})}]{MR3441527}
\bibinfo{author}{Betsakos, D.} (\bibinfo{year}{2016}{\natexlab{b}}).
\newblock \bibinfo{title}{On the asymptotic behavior of the trajectories of
  semigroups of holomorphic functions}.
\newblock {\it \bibinfo{journal}{J. Geom. Anal.}\/},  {\it
  \bibinfo{volume}{26}\/}, \bibinfo{pages}{557--569}.
  \DOIprefix\doi{10.1007/s12220-015-9562-1}.
\bibitem[{Bracci et~al.(2018{\natexlab{a}})Bracci, Contreras,
  D{\'\i}az-Madrigal \& Gaussier}]{bracci2018non}
\bibinfo{author}{Bracci, F.}, \bibinfo{author}{Contreras, M.~D.},
  \bibinfo{author}{D{\'\i}az-Madrigal, S.}, \& \bibinfo{author}{Gaussier, H.}
  (\bibinfo{year}{2018}{\natexlab{a}}).
\newblock \bibinfo{title}{Non-tangential limits and the slope of trajectories
  of holomorphic semigroups of the unit disc}.
\newblock {\it \bibinfo{journal}{arXiv preprint arXiv:1804.05553v2}\/}, .
\bibitem[{Bracci et~al.(2018{\natexlab{b}})Bracci, Contreras,
  D{\'\i}az-Madrigal, Gaussier \& Zimmer}]{bracci2018asy}
\bibinfo{author}{Bracci, F.}, \bibinfo{author}{Contreras, M.~D.},
  \bibinfo{author}{D{\'\i}az-Madrigal, S.}, \bibinfo{author}{Gaussier, H.}, \&
  \bibinfo{author}{Zimmer, A.} (\bibinfo{year}{2018}{\natexlab{b}}).
\newblock \bibinfo{title}{Asymptotic behavior of orbits of holomorphic
  semigroups}.
\newblock {\it \bibinfo{journal}{arXiv preprint arXiv:1810.07947}\/}, .
\bibitem[{Carath\'{e}odory(1954)}]{MR0060009}
\bibinfo{author}{Carath\'{e}odory, C.} (\bibinfo{year}{1954}).
\newblock {\it \bibinfo{title}{Theory of functions of a complex variable.
  {V}ol. 1}\/}.
\newblock \bibinfo{publisher}{Chelsea Publishing Co., New York, N. Y.}
\newblock \bibinfo{note}{Translated by F. Steinhardt}.
\bibitem[{Contreras \& D\'{i}az-Madrigal(2005)}]{MR2225072}
\bibinfo{author}{Contreras, M.~D.}, \& \bibinfo{author}{D\'{i}az-Madrigal, S.}
  (\bibinfo{year}{2005}).
\newblock \bibinfo{title}{Analytic flows on the unit disk: angular derivatives
  and boundary fixed points}.
\newblock {\it \bibinfo{journal}{Pacific J. Math.}\/},  {\it
  \bibinfo{volume}{222}\/}, \bibinfo{pages}{253--286}.
  \DOIprefix\doi{10.2140/pjm.2005.222.253}.
\bibitem[{Contreras et~al.(2015)Contreras, D\'{i}az-Madrigal \&
  Gumenyuk}]{MR3318309}
\bibinfo{author}{Contreras, M.~D.}, \bibinfo{author}{D\'{i}az-Madrigal, S.}, \&
  \bibinfo{author}{Gumenyuk, P.} (\bibinfo{year}{2015}).
\newblock \bibinfo{title}{Slope problem for trajectories of holomorphic
  semigroups in the unit disk}.
\newblock {\it \bibinfo{journal}{Comput. Methods Funct. Theory}\/},  {\it
  \bibinfo{volume}{15}\/}, \bibinfo{pages}{117--124}.
  \DOIprefix\doi{10.1007/s40315-014-0092-9}.
\bibitem[{Elin \& Jacobzon(2014)}]{MR3240989}
\bibinfo{author}{Elin, M.}, \& \bibinfo{author}{Jacobzon, F.}
  (\bibinfo{year}{2014}).
\newblock \bibinfo{title}{Parabolic type semigroups: asymptotics and order of
  contact}.
\newblock {\it \bibinfo{journal}{Anal. Math. Phys.}\/},  {\it
  \bibinfo{volume}{4}\/}, \bibinfo{pages}{157--185}.
  \DOIprefix\doi{10.1007/s13324-014-0084-y}.
\bibitem[{Elin et~al.(2010)Elin, Khavinson, Reich \& Shoikhet}]{MR2731701}
\bibinfo{author}{Elin, M.}, \bibinfo{author}{Khavinson, D.},
  \bibinfo{author}{Reich, S.}, \& \bibinfo{author}{Shoikhet, D.}
  (\bibinfo{year}{2010}).
\newblock \bibinfo{title}{Linearization models for parabolic dynamical systems
  via {A}bel's functional equation}.
\newblock {\it \bibinfo{journal}{Ann. Acad. Sci. Fenn. Math.}\/},  {\it
  \bibinfo{volume}{35}\/}, \bibinfo{pages}{439--472}.
  \DOIprefix\doi{10.5186/aasfm.2010.3528}.
\bibitem[{Elin \& Shoikhet(2010)}]{MR2683159}
\bibinfo{author}{Elin, M.}, \& \bibinfo{author}{Shoikhet, D.}
  (\bibinfo{year}{2010}).
\newblock {\it \bibinfo{title}{Linearization models for complex dynamical
  systems}\/} volume \bibinfo{volume}{208} of {\it \bibinfo{series}{Operator
  Theory: Advances and Applications}\/}.
\newblock \bibinfo{publisher}{Birkh\"{a}user Verlag, Basel}.
\newblock \DOIprefix\doi{10.1007/978-3-0346-0509-0}.
\bibitem[{Garnett \& Marshall(2008)}]{MR2450237}
\bibinfo{author}{Garnett, J.~B.}, \& \bibinfo{author}{Marshall, D.~E.}
  (\bibinfo{year}{2008}).
\newblock {\it \bibinfo{title}{Harmonic measure}\/} volume~\bibinfo{volume}{2}
  of {\it \bibinfo{series}{New Mathematical Monographs}\/}.
\newblock \bibinfo{publisher}{Cambridge University Press, Cambridge}.
\newblock \bibinfo{note}{Reprint of the 2005 original}.
\bibitem[{Pommerenke(1992)}]{MR1217706}
\bibinfo{author}{Pommerenke, C.} (\bibinfo{year}{1992}).
\newblock {\it \bibinfo{title}{Boundary behaviour of conformal maps}\/} volume
  \bibinfo{volume}{299} of {\it \bibinfo{series}{Grundlehren der Mathematischen
  Wissenschaften [Fundamental Principles of Mathematical Sciences]}\/}.
\newblock \bibinfo{publisher}{Springer-Verlag, Berlin}.
\newblock \DOIprefix\doi{10.1007/978-3-662-02770-7}.
\bibitem[{Ransford(1995)}]{MR1334766}
\bibinfo{author}{Ransford, T.} (\bibinfo{year}{1995}).
\newblock {\it \bibinfo{title}{Potential theory in the complex plane}\/}
  volume~\bibinfo{volume}{28} of {\it \bibinfo{series}{London Mathematical
  Society Student Texts}\/}.
\newblock \bibinfo{publisher}{Cambridge University Press, Cambridge}.
\newblock \DOIprefix\doi{10.1017/CBO9780511623776}.
\bibitem[{Shoikhet(2001)}]{MR1849612}
\bibinfo{author}{Shoikhet, D.} (\bibinfo{year}{2001}).
\newblock {\it \bibinfo{title}{Semigroups in geometrical function theory}\/}.
\newblock \bibinfo{publisher}{Kluwer Academic Publishers, Dordrecht}.
\newblock \DOIprefix\doi{10.1007/978-94-015-9632-9}.

\end{thebibliography}
  
\end{document}